\documentclass[12pt,a4paper]{article}
\usepackage[T2A]{fontenc}
\usepackage[cp1251]{inputenc}
\usepackage[russian]{babel}
\usepackage{srcltx}\usepackage{graphicx,color}
\usepackage[all]{xy}
\usepackage{latexsym,amsfonts,amssymb,amsmath,longtable,%amsthm,
mathrsfs,theorem,cite
}
\renewcommand{\geq}{\geqslant}

\newcommand{\qed}{\hfill{$\Box$}}
\theoremstyle{plain}
{%\theorembodyfont{\itshape}
\newtheorem{Lemma}{{\bfseries Лемма}}}
{\theorembodyfont{\itshape}

}

\theorembodyfont{\rm}

\DeclareMathOperator{\Syl}{Syl} 
\DeclareMathOperator{\prn}{\,prn\,}\DeclareMathOperator{\PSp}{PSp}

\title{\vspace{-1cm} \hfill{%\normalsize УДК 512.542}{
\fontfamily{cmr} \fontseries{bx} \selectfont \\ \vspace{1cm} On the pronormality of subgroups of odd indices in finite simple symplectic groups \\ {\small (in Russian, English abstract)}}
\thanks{The work is supported by the grant of the President of the Russian Federation for young scientists (grant no. MK-6118.2016.1), by the Integrated Program for Fundamental Research of the Ural Branch of the Russian Academy of Sciences (project no. 15-16-1-5), by the Program for State Support of Leading Universities of the Russian Federation (agreement
no. 02.A03.21.0006 of August 27, 2013). The second author is a winner of the competition of young mathematicians of the Dmitry Zimin Foundation ``Dynasty'' (in 2013 year). The third author is supported by Chinese Academy of Sciences President's International Fellowship Initiative (PIFI), grant no.~2016VMA078.}}
\date{}
\author{\bf   Anatoly~S.~Kondrat'ev, Natalia~V.~Maslova, Danila~O.~Revin}

\begin{document}
%%%%%%%%%%%%%%%%%%%%%%%%%%%%%%%%%%%%%%%%%%%%%%%%%%%%%%%%%%%%%%%%%%%%%%%%%%%%%%%%%%%

%\renewcommand{\qed}{}
%\renewcommand{\proofname}{{ {\rm ДОКАЗАТЕЛЬСТВО}}}

%\thispagestyle{empty}
%\vspace{-2cm}

\maketitle
\pagenumbering{arabic}
\begin{center}
{\bf {\small Abstract}}
\end{center}

{\small A subgroup $H$ of a group $G$ is said to be {\it pronormal} in $G$ if $H$ and $H^g$ are conjugate in $\langle H, H^g \rangle$ for every $g \in G$. In [Sib. Math. J. 2012. Vol. 53, no.~3], the following conjecture was formulated by Evgeny~P.~Vdovin and the third author.

\medskip
{\bf Conjecture.} {\sl All subgroups of odd indices are pronormal in all finite simple groups.}
\medskip

The conjecture was verified by authors for many families of finite simple groups in [Sib. Math. J. 2015.  Vol. 56, no.~6]. Namely, the following theorem was proved.

\medskip

{\bf Theorem.}
{\sl All subgroups of odd indices are pronormal in the following finite simple groups{\rm:}
$A_n$, where $n\ge 5${\rm;}
sporadic groups{\rm;}
groups of Lie type over fields of characteristic $2${\rm;}
$L_{2^n}(q)${\rm;}
$U_{2^n}(q)${\rm;}
$S_{2n}(q)$, where  $q\not\equiv\pm3 \pmod 8${\rm;}
$O_n(q)${\rm;}
exeptional groups of Lie type not isomorphic to $E_6(q)$ or ${}^2E_6(q)$.}

\medskip

In [Proc. Steklov Inst. Math., to appear, Theorem~1] authors proved that, if a group $G$ has a normal abelian subgroup $V$ and a subgroup $H$ such that $G=HV$, then $H$ is pronormal in $G$ if and only if $U=N_U(H)[H,U]$ for any $H$-invariant subgroup $U$ of the group $V$. Using this fact, in [Proc. Steklov Inst. Math., to appear, Theorem~2] it was proved that Conjecture fails. Precisely, a finite simple symplectic group $PSp_{6n}(q)$ with $q \equiv \pm 3 \pmod 8$ contains a nonpronormal subgroup of odd index.
In view of the above results, the following problem naturally arises.

\medskip	
{\bf Problem.} Classify finite simple groups in which all subgroups of odd indices are pronormal.

\medskip

Using [Proc. Steklov Inst. Math., to appear, Theorem~2] we prove the following theorem.
\medskip

 {\bf Theorem 1.} Let $G=PSp_{2n}(q)$, where $q \equiv \pm 3 \pmod 8$ and $n \not \in \{2^m, 2^m(2^{2k}+1) \mid m, k \in \mathbb{N} \cup \{0\}\}$. Then $G$ contains a nonpronormal subgroup of odd index.
\medskip

The main result of this paper is the following theorem.
\medskip

{\bf Theorem 2.} Let $G=PSp_{2^n}(q)$, where $n \ge 2$ and $q \equiv \pm 3 \pmod 8$. Then any subgroup of odd index of $G$ is pronormal in $G$.

}



\begin{center}
{\bf {\small Аннотация}}
\end{center}

{\small Подгруппа $H$ группы $G$ называется пронормальной, если для любого элемента $g\in G$ подгруппы $H$ и $H^g$ сопряжены в
подгруппе $\langle H, H^g\rangle$.

В работе Е.П.Вдовина и третьего автора [Сиб. матем. журн. 2012. Т~53, \No~3] была высказана гипотеза  о том, что
подгруппа нечетного индекса в конечной простой группе всегда пронормальна. Недавно [Сиб. матем. журн. 2015. Т~56, \No~6] авторы подтвердили эту гипотезу
 для всех конечных простых групп, за исключением $PSL_n(q)$, $PSU_n(q)$, $E_6(q)$ и ${}^2E_6(q)$, где $q$ во всех случаях нечетно и $n$ не является степенью числа $2$ в первых двух случаях, а также $PSp_{2n}(q)$, где $q \equiv \pm 3 \pmod 8$.

Однако в работе [Тр. Ин-та математики и механики УрО РАН. 2016. Т.~22, \No~1] авторами было доказано,  что при $q\equiv\pm 3\pmod 8$ и $n\equiv 0\pmod 3$ простая симплектическая группа $PSp_{2n}(q)$ содержит непронормальную
подгруппу нечетного индекса. Тем самым гипотеза  о пронормальности подгруппы нечетного индекса в конечной простой группе была
опровергнута.

Как естественное расширение данной гипотезы возникает проблема классификации конечных неабелевых простых групп, в которых любая подгруппа нечетного индекса пронормальна. В настоящей работе мы продолжаем изучать эту проблему для симплектической простой группы~$PSp_{2n}(q)$ при $q \equiv \pm 3 \pmod 8$ (в~отсутствие этого ограничения подгруппы нечетных индексов пронормальны). Доказывается, что если  $n$ не является числом вида $2^m$ или $2^m(2^{2k}+1)$,
то данная группа содержит непронормальную подгруппу нечетного индекса. В то же время в статье доказано, что если $n=2^m$, то все подгруппы нечетных индексов в группе~$PSp_{2n}(q)$ пронормальны. Для случая $n=2^m(2^{2k}+1)$ и $q \equiv \pm 3 \pmod 8$ вопрос о пронормальности подгрупп нечетных индексов в группе $PSp_{2n}(q)$ пока остается открытым.

}

%%%%%%%%%%%%%%%%%%%%%%%%%%%%%%%%%%%%%%%%%%%%%%%%%%%%%%%%%%%%%%%%%%%%%%%

\section*{Введение}

%%%%%%%%%%%%%%%%%%%%%%%%%%%%%%%%%%%%%%%%%%%%%%%%%%%%%%%%%%%%%%%%%%%%%%%%

В соответствии с определением Ф.\,Холла, подгруппа $H$ группы $G$ называется {\em пронормальной},
если для любого элемента $g\in G$ подгруппы $H$~и~$H^g$ сопряжены в подгруппе $\langle H, H^g\rangle$.

В дальнейшем мы рассматриваем только конечные группы, и в связи с этим термин <<группа>> употребляется нами в значении <<конечная группа>>.

Подгруппа $H$ конечной группы $G$ называется {\it холловой}, если ее порядок $|H|$ и индекс $|G:H|$ взаимно просты.

\medskip

В работе \cite{Pronorm_Hall_Simple} было доказано, что холловы подгруппы в конечных простых группах пронормальны, и на основании анализа доказательства была высказана следующая

\medskip

{\bf Гипотеза 1}{{\rm\cite[гипотеза 1]{Pronorm_Hall_Simple}}}. В конечных простых группах подгруппы нечетных индексов пронормальны.

\medskip

Эта гипотеза была подтверждена авторами \cite[теорема]{KMR_1} для конечных простых групп, за исключением $PSL_n(q)$, $PSU_n(q)$, $E_6(q)$ и ${}^2E_6(q)$, где $q$ во всех случаях нечетно и $n$ не является степенью числа $2$ в первых двух случаях, а также $PSp_{2n}(q)$, где $q \equiv \pm 3 \pmod 8$.

Однако с помощью полученного в \cite[теорема~1]{KMR_2} критерия пронормальности добавлений к абелевым нормальным подгруппам конечных групп авторам удалось построить примеры непронормальных подгрупп нечетных индексов в группах $PSp_{6n}(q)$, где $q \equiv \pm 3 \pmod 8$ (см. \cite[теорема~2]{KMR_2}). Таким образом, гипотеза~1 была опровергнута.

Как естественное расширение гипотезы 1 в свете полученных результатов возникает

\medskip
{\bf Проблема 1.} Классифицировать неабелевы простые группы, в которых подгруппы нечетных индексов пронормальны.
\medskip

В настоящей работе мы продолжаем исследование проблемы~1 для случая симплектических групп.

На самом деле, результат \cite[теорема~2]{KMR_2} позволяет исследовать проблему 1 для значительно б\'{о}льшего массива симплектических групп, чем это было сделано в \cite{KMR_2}. Мы докажем следующее утверждение.

\medskip
{\bf Теорема 1.} {\it Пусть $G=PSp_{2n}(q)$, где $q \equiv \pm 3 \pmod 8$ и $n$ не является числом вида $2^w$ или $2^w({2^{2k}+1})$. Тогда $G$ содержит непронормальную подгруппу нечетного индекса.}
\medskip

Кроме того, мы полностью исследуем один из двух оставшихся случаев симплектических групп, которые не накрываются теоремой~1, а именно случай групп $PSp_{2n}(q)$, где  $q \equiv \pm 3 \pmod 8$ и $n=2^w$. Основным результатом настоящей работы является

\medskip
{\bf Теорема 2.} {\it Подгруппы нечетных индексов пронормальны в группах $PSp_{2n}(q)$, где $n=2^w \ge 2$ и $q \equiv \pm 3 \pmod 8$. }
\medskip

Из теоремы 2 и основного результата работы \cite{KMR_1} получаем

\medskip
{\bf Следствие.} {\it Подгруппы нечетных индексов пронормальны в простой классической группе, у которой размерность или порядок основного поля является степенью числа $2$. }
\medskip

{\bf Замечание.} Проблема 1 остается открытой для групп $PSp_{2n}(q)$, где $q \equiv \pm 3 \pmod 8$ и $n=2^m(2^{2k}+1)$, а также для групп $E_6(q)$, ${}^2E_6(q)$, $PSL_n(q)$ и $PSU_n(q)$, где во всех случаях $q$ нечетно и $n$ не является степенью числа $2$.

\section*{Вспомогательные результаты}

Наши терминология и обозначения, в основном, стандартны, их можно найти в \cite{Atlas,KL}.
Запись $H\prn G$ используется для сокращения ``$H$ является пронормальной подгруппой в группе $G$''.

Для группы $G$ и подмножества $\pi$ множества всех простых чисел через $Soc(G)$, $O_\pi(G)$ и $Z(G)$ обозначаются цоколь (подгруппа, порожденная всеми минимальными неединичными нормальными подгруппами), $\pi$-радикал (наибольшая нормальная $\pi$-подгруппа) и центр группы $G$ соответственно. Через
${\rm Syl}_p(G)$ обозначается множество силовских $p$-подгрупп группы~$G$. Как обычно, через $\pi'$ обозначается множество всех простых чисел, не принадлежащих~$\pi$. Кроме того, если $n$~--- натуральное число, то $n_\pi$~--- наибольший натуральный делитель числа $n$ такой, что все простые делители $n_\pi$ принадлежат~$\pi$. Для конечной группы $G$ будем, как это обычно принято, писать $O(G)$ вместо $O_{2'}(G)$.

Зафиксируем некоторые обозначения, связанные с классическими группами. Пусть $q$ --- натуральная степень простого числа и $G$ --- одна из конечных простых классических групп $PSL_n(q)$, $PSU_n(q)$, $PSp_n(q)$ для четного~$n$, $P\Omega_n(q)$ для нечетных $n$ и $q$ и $P\Omega_n^\varepsilon(q)$ для четного~$n$, где $\varepsilon \in \{+,-\}$. Будем обозначать через $V$ векторное пространство размерности $n$ над полем~$F$ с соответствующей билинейной или квадратичной формой, ассоциированное с группой $G$, где $F=\mathbb{F}_q$ для линейных, симплектических и ортогональных групп и $F=\mathbb{F}_{q^2}$ для унитарных групп.

Здесь и далее мы используем следующие обозначения \cite{KL}:  $PSL_n(q)=L_n(q)=L_n^+(q)$, $PSU_n(q)=U_n(q)=L_n^-(q)$, $PSp_n(q)=S_n(q)$,  $E^{+}_6(q)=E_6(q)$ и  $E^{-}_6(q)={}^2E_6(q)$.

\begin{Lemma}{\rm \cite[лемма 5]{Pronorm_Hall_Simple}} \label{NormSyl} Пусть $G$ --- конечная группа и $H\le G$.
Предположим также, что подгруппа $H$ содержит некоторую силовскую подгруппу $S$ группы $G$. Тогда следующие утверждения эквивалентны:
\begin{itemize}
\item[$(1)$] $H\prn G$;
\item[$(2)$] подгруппы $H$ и $H^g$ сопряжены в $\langle H, H^g\rangle$ для любого $g\in N_G(S)$.
\end{itemize}
\end{Lemma}

Таким образом, для доказательства пронормальности подгруппы $H$ нечетного индекса в
группе $G$ достаточно установить сопряженность подгрупп $H$ и $H^g$  в $\langle H, H^g\rangle$ для любого нетривиального
элемента
$g\in N_G(S)$ нечетного порядка, где $S$ --- фиксированная силовская $2$-подгруппа группы $H$ (и, следовательно, группы $G$).

\begin{Lemma}\label{N_G(S)=S} Пусть $p$~--- простое число и $\mathfrak{X}_p$~--- класс всех конечных групп, в которых силовская $p$-подгруппа совпадает со своим нормализатором. Тогда
 \begin{itemize}
  \item[$(1)$] если $X\unlhd Y$, причем $X\in \mathfrak{X}_p$ и $Y/X\in \mathfrak{X}_p$, то $Y\in \mathfrak{X}_p${\rm;}
  \item[$(2)$] если $X\le Y\in \mathfrak{X}_p$ и индекс $|Y:X|$ не делится на $p$, то $X\prn Y$.
  \end{itemize}
\end{Lemma}

\noindent{\sc Доказательство.} Докажем утверждение  (1). Пусть  $X\unlhd Y$, причем $Y/X\in \mathfrak{X}_p$ и $X\in \mathfrak{X}_p$. Пусть также $S\in \Syl_p(Y)$. Справедливы равенства $$N_Y(SX)=SX\text{ и }N_X(S\cap X)=S\cap X,$$ поскольку $SX/X\in \Syl_p(Y/X)$ и $S\cap X\in \Syl_p(X)$. Покажем, что  $N_Y(S)= S$. Это так, поскольку $N_Y(S)\le N_Y(SX)=SX$ и, следовательно, $$N_Y(S)=N_{SX}(S)=SN_X(S)\,\,\,\text{ и }\,\,\,N_X(S)\le N_X(S\cap X)=S\cap X.$$

Утверждение (2) следует из леммы~\ref{NormSyl} и того факта, что подгруппа индекса, не кратного~$p$, в группе $Y\in \mathfrak{X}_p$ содержит некоторую силовскую $p$-подгруппу $S$ группы~$Y$, а значит и ее нормализатор в~$Y$. \qed\medskip

\begin{Lemma}{\rm \cite[следствие]{Kond1}} \label{Norm_Syl_Simple}  Пусть $G$ --- конечная простая неабелева группа и $S\in\Syl_2(G)$. Тогда $N_G(S)=S$, за исключением следующих случаев:
\begin{itemize}

\item[$(1)$]  $G\cong  J_2,~J_3,~Suz$ или $F_5$ и $|N_G(S):S|=3$;

\item[$(2)$]  $G\cong  {^2}G_2(3^{2n+1})$ или $J_1$ и $N_G(S)\cong 2^3.7.3<{\rm Hol}(2^3)$;

\item[$(3)$] $G$ --- группа лиева типа над полем характеристики $2$ и $N_G(S)$ --- подгруппа Бореля в группе $G$.

\item[$(4)$] $G\cong L_2(q)$, где $3<q \equiv  \pm 3\pmod 8$ и $N_G(S)\cong A_4$.

\item[$(5)$]  $G \cong E^{\eta}_6(q)$, $\eta=\pm$, $q$ нечетно и
$N_G(S)=S\times C$, где $C$ --- неединичная циклическая группа порядка
$(q-\eta 1)_{2'}/(q-\eta 1,3)_{2'}$.

\item[$(6)$]  $G\cong S_{2n}(q)$, $n\geq 2$, $q\equiv\pm 3\pmod 8$ и
факторгруппа $N_G(S)/S$ изоморфна
элементарной абелевой $3$-группе порядка $3^t$, где число $t$ находится из
двоичного разложения $$n=2^{s_1}+\dots+2^{s_t},~s_1>\dots>s_t\geq 0.$$

\item[$(7)$] $G\cong L_n^\eta(q)$, $n\geq 3$, $\eta=\pm$, $q$ нечетно и $$N_G(S)\cong S\times C_1\times \dots\times C_{t-1},$$ где число $t$ находится из двоичного разложения
    $$n=2^{s_1}+\dots+2^{s_t},
        $$
    $s_1>\dots>s_t>0,$ и $C_1,\dots C_{t-2},C_{t-1}$ --- циклические группы  порядков \\
${(q-\eta 1)}_{2'},\dots,{(q-\eta 1)}_{2'},{(q-\eta 1)}_{2'}/{(q-\eta 1,n)} _{2'}$ соответственно.
\end{itemize}
\end{Lemma}

%\begin{Lemma}{\rm \cite[теорема II.8.27]{Hupp}} \label{Hupp}
%Пусть $H$ --- подгруппа группы $PSL_2(p^s)$, где $p$~--- простое число и $s$~--- положительное целое число.
%Обозначим $(p^s-1,2)$ через~$d$.
%Тогда $H$ изоморфна одной из следующих групп:
%\begin{itemize}%
	%\item[$(1)$] элементарной абелевой $p$-группе порядка, не превосходящего $p^s$;
	
	%\item[$(2)$] циклической группе порядка $z$, где $z$ --- делитель числа $\frac{p^s \pm 1}{d}$;
	
	%\item[$(3)$] группе диэдра порядка $2z$,  где число $z$ --- делитель числа $\frac{p^s \pm 1}{d}$;
	
	%\item[$(4)$] $A_4$, при этом либо $p$ нечётно, либо $p=2$ и $s$ чётно;
	
	%\item[$(5)$] $S_4$, при этом $p^{2s}-1$ делится на $16$;
	
	%\item[$(6)$] $A_5$, при этом $p=5$ или $p^{2s}-1$ делится на $5$;
	
	%\item[$(7)$] полупрямому произведению элементарной абелевой группы порядка $p^m$ и циклической группы порядка $t$ такого, что $t$ делит $p^{(m,s)} - 1$;
	
	%\item[$(8)$] группе $PSL_2(p^m)$, если $m$ делит $s$;

%\item[$(9)$] группе $PGL_2(p^m)$, если $2m$ делит $s$.
%\end{itemize}
%\end{Lemma}

\begin{Lemma}{\rm \cite[лемма~1]{Maz}} \label{Maz} Если группа Фробениуса $FC$ с ядром $F$ и циклическим дополнением $C=\langle g\rangle$ порядка $n$ точно действует на векторном пространстве $V$ над полем ненулевой характеристики, взаимно простой с $|F|$, то минимальный многочлен элемента $g$ как линейного преобразования пространства $V$ равен $\lambda^n-1$.
\end{Lemma}

\begin{Lemma}\label{OddOrderSubgroupsPSL}
Пусть $H$ --- подгруппа нечетного порядка группы~$G$, где $$SL_2(q)\le G\le GL_2(q)\text{ или } PSL_2(q)\le G \le PGL_2(q),$$ причем $(|H|,q)=1$, и $g\in N_G(H)$~--- элемент нечетного порядка. Тогда
\begin{itemize}
  \item[$(1)$] группа $H$ абелева;
  \item[$(2)$]  $g\in C_G(H)$.
\end{itemize}
\end{Lemma}

\noindent{\sc Доказательство.} Для группы $G=PSL_2(q)$ данное утверждение легко вывести из списка максимальных подгрупп группы $G$  \cite[теорема II.8.27]{Hupp}. Мы приведем однако здесь другое доказательство, охватывающее все остальные случаи и основанное на элементах теории представлений и линейной алгебры.

Допустим, $PSL_2(q)\le G \le PGL_2(q)$. Так как  подгруппы $H$ и  $\langle H,g\rangle$ имеют нечетные порядки, эти подгруппы содержатся в~$PSL_2(q)$. Более того, в $SL_2(q)\le GL_2(q)$ у этих подгрупп найдутся изоморфные им прообразы. Следовательно, нам достаточно доказать лемму для случая, когда $G=GL_2(q)$.

Докажем (1). По теореме Машке \cite[теорема~(1.9)]{Isaacs} группа $H$ вполне приводима. Следовательно, если $H$ приводима, она сопряжена с подгруппой группы диагональных матриц, и утверждение (1) верно. Поэтому считаем, что подгруппа $H$ неприводима. Более того, заменяя в случае необходимости поле $\mathbb{F}_q$ на поле разложения  $\mathbb{F}_{q^n}$ группы $H$ и пользуясь вложением $GL_2(q)\hookrightarrow GL_2(q^n)$, мы, не уменьшая общности, можем считать тождественное вложение $H$ в $G$ точным абсолютно неприводимым представлением группы $H$. По~\cite[теорема~(15.13)]{Isaacs} степень этого представления, равная~2, совпадает со степенью некоторого неприводимого комплексного представления, и, следовательно, по~\cite[теорема~(3.11)]{Isaacs} делит нечетное число~$|H|$. Противоречие.

Докажем (2). Пусть характеристика поля  $\mathbb{F}_q$ равна $p$. Элемент $g$ можно представить в виде $g=xy$, где $x,y\in\langle g\rangle$, причем $x$~--- $p$-элемент, а $y$~--- $p'$-элемент. Теперь, заменяя $H$ на подгруппу $\langle H,y\rangle$ (порядкок которой равен $|H||y|$ и которая, ввиду утверждения (1), абелева), а $g$ на $x$, мы можем считать, что $g$~--- $p$-элемент. Более того, поскольку силовская $p$-подгруппа группы $G$ изоморфна аддитивной группе поля  $\mathbb{F}_q$, мы можем считать, что~$|g|=p$.  Ввиду~\cite[теорема~4.34]{Isaacs1} и абелевости подгруппы~$H$, имеем $H=F\times K$, где $F=[H,\langle g\rangle]$ и $K=C_H(g)$. При этом $C_F(g)=F\cap K=1$, поэтому группа $\langle F, g\rangle=F\langle g\rangle$ является группой Фробениуса с ядром $F$ и циклическим дополнением $\langle g\rangle$ порядка~$p$.   Применяя лемму \ref{Maz}, получаем противоречие, так как минимальный многочлен $\lambda^p-1$ элемента $g$  делит его характеристический многочлен $\det(g-\lambda I)$ степени~2, а $p=|g|$ нечетно. \qed\medskip

\begin{Lemma}{\rm \cite[лемма~3]{KMR_1}} \label{Quot} Пусть $H$ --- подгруппа, $N$~--- нормальная подгруппа группы $G$ и $\overline{\phantom{g}}:G\rightarrow G/N$~---
естественный  эпиморфизм. Справедливы следующие утверждения:
\begin{itemize}
\item[$(1)$] если $H\prn G$, то $\overline{H}\prn\overline{G}$;
\item[$(2)$] если $N\le H$ и  $\overline{H}\prn\overline{G}$, то $H\prn G$.
\end{itemize}
 В частности, подгруппа $H$ нечетного индекса пронормальна в $G$ тогда и только тогда, когда подгруппа $H/O_2(G)$
пронормальна в $G/O_2(G)$.
\end{Lemma}

%\begin{Lemma}{\rm \cite[лемма 4]{Pronorm_Hall_Simple}} \label{Instead} Пусть $H$ --- подгруппа  группы $G$ и $g\in G$. Допустим, что для некоторого элемента
%$y\in \langle H, H^g\rangle$ подгруппы $H^y$ и $H^g$ сопряжены в $ \langle H^y, H^g\rangle$. Тогда подгруппы  $H$ и $H^g$ сопряжены в $ \langle H, H^g\rangle$.
%\end{Lemma}

\begin{Lemma}{\rm \cite[лемма~5]{KMR_1}} \label{Overgroup} Пусть $H$ и $M$ --- подгруппы  группы $G$, причем $H\le M$. Справедливы следующие утверждения:
\begin{itemize}
\item[$(1)$] если $H\prn G$, то $H\prn M$;
\item[$(2)$] если $S\le H$ для некоторой силовской подгруппы $S$ группы $G$, причем  ${N_G(S)}\le{M}$ и $H\prn M$, то $H\prn G$.
\end{itemize}\end{Lemma}

%\begin{Lemma} \label{AS} Пусть $\pi$~--- некоторое множество простых чисел. Допустим, что ${A\unlhd G}$ и $G/A$~--- $\pi$-группа.
%Предположим также, что индекс некоторой подгруппы $H$ группы $G$ не делится на числа из $\pi$ и $(H\cap A)\prn A$. Тогда  $H\prn G$.
%\end{Lemma}

\begin{Lemma}{\rm \cite[лемма~3]{KMR_2}} \label{TransitiveAction} Пусть $H$~--- транзитивная группа подстановок степени $n$ и $A$~--- группа. Определим действие группы $H$ на группе $V=A^n$ по правилу
$$
(x_1,\dots, x_n)^\pi=(x_{1\pi^{-1}},\dots, x_{n\pi^{-1}}) \text{ для любых } (x_1,\dots, x_n)\in V \text{ и } \pi\in H.
$$
Тогда $$C_V(H)=\{(x_1,\dots, x_n)\in V\mid x_1=\dots= x_n\}\cong A.$$
\end{Lemma}

\begin{Lemma}{\rm \cite[теорема~1]{KMR_2}} \label{CrterionOfProSub} Пусть $H$ и $V$~---  подгруппы группы $G$ такие, что $V$~--- абелева нормальная подгруппа в $G$ и
 $G=HV$. Тогда следующие утверждения равносильны{\rm:}
\begin{itemize}
\item[$(1)$] $H \prn G${\rm;}
\item[$(2)$]  $U=N_U(H)[H,U]$ для любой $H$-инвариантной подгруппы $U\le V$.
\end{itemize}
\end{Lemma}

\begin{Lemma}{\rm \cite[теорема~2]{KMR_2}} \label{NonProNorm}  Группа $\PSp_{6n}(q)$ при $q\equiv\pm 3\pmod 8$  содержит непронормальную подгруппу нечетного индекса.
\end{Lemma}

\begin{Lemma}{\rm \cite[лемма~8]{KMR_2}} \label{AwrH} Пусть $H$~--- группа подстановок степени $n$ и $A$~--- конечная группа. Определим действие группы $H$ на группе $V=A^n$ по правилу
$$
(x_1,\dots, x_n)^\pi=(x_{1\pi^{-1}},\dots, x_{n\pi^{-1}}) \text{ для любых } (x_1,\dots, x_n)\in V \text{ и } \pi\in H.
$$
Предположим, что $H$ содержит транзитивную подгруппу $K$ такую, что $$(|A|,|K|)=1.$$ Тогда для любой $H$-инвариантной подгруппы $U$ группы $V$ выполнено равенство
$$
U=C_U(H)[H,U].
$$
В частности, $H \prn VH$.
\end{Lemma}

\begin{Lemma}\label{DwrSn} Пусть выполнено одно из следующих утверждений:
\begin{itemize}
  \item[$($а$)$] $G=S_n$~--- симметрическая группа степени $n${\rm;}
  \item[$($б$)$] $G$~--- диэдральная группа;
  \item[$($в$)$] $G=D\wr S_n$~--- сплетение диэдральной и симметрической групп.\end{itemize}
Тогда{\rm:}
\begin{itemize}
  \item[$(1)$] силовская $2$-подгруппа группы $G$ совпадает со своим нормализатором{\rm;}
  \item[$(2)$] любая подгруппа нечетного индекса пронормальна в $G$. \end{itemize}
\end{Lemma}

\noindent{\sc Доказательство.} Достаточно доказать утверждение (1) в случаях (а) и (б). Тогда справедливость утверждения (1) в случае (в)  и утверждение (2) будут следовать из леммы~\ref{N_G(S)=S}. Для случая (а) см.~\cite[лемма 4]{CF}. Рассмотрим случай (б). Путь $G$~--- диэдральная группа, $S\in\Syl_2(G)$ и $U=O(G)$. Тогда $G=SU$, поэтому $$N_G(S)=SN_U(S)=SC_U(S).$$ Покажем, что $C_U(S)=1$ и тем самым докажем утверждение (1). Из определения диэдральной группы следует, что группа $U$ циклическая и $S$ содержит инволюцию $t$ такую, что $u^t=u^{-1}$ для всех $u\in U$. Допустим, $1\ne u\in C_U(S)$. Тогда $u$~--- нетривиальный элемент нечетного порядка и $$u\ne u^{-1}=u^t=u.$$ Получаем противоречие. \qed\medskip

\begin{Lemma}\label{3wrSn} Пусть $G=A\wr B=VB$~--- естественное подстановочное сплетение циклической группы $A$ порядка $3$ и симметрической группы $B=S_n$, где $n=2^w$, через $V=V_1\times \ldots \times V_n$, где $V_i \cong A$, обозначена база сплетения, и $H$ --- подгруппа нечетного индекса в $G$. Тогда $H \prn G$.
\end{Lemma}

\noindent{\sc Доказательство.} Пусть $S \in \Syl_2(G)$ и $S \le H$. По теореме Силова можно считать, что $S \le B$.

%Сначала покажем, что $N_G(S) \le HV$. Заметим, что подгруппа $VS$ нормальна в $N_G(S)V$, т.~е. $N_G(S)V \le N_G(SV)$. Пусть $\bar{ } :G \rightarrow G/V$ --- естественный гомоморфизм. Заметим, что $\overline{N_G(SV)} = N_{\overline{G}}(\overline{S})$. Поскольку порядок $V$ нечетен, из того, что $S \in \Syl_2(H)$, следует, что $\overline{S} \in \Syl_2(\overline{H})$. Из \cite[лемма~4]{CF} следует, что $N_{\overline{G}}(\overline{S})=\overline{S} \le \overline{H}$, поэтому $N_G(SV) \le HV$. Отсюда $$N_G(S) \le N_G(S)V \le N_G(SV) \le HV.$$

Из леммы~\ref{DwrSn}(1) следует, что $N_G(SV)=SV \le HV$. Поэтому $$N_G(S)\le N_G(SV) \le HV.$$

Теперь ввиду леммы \ref{Overgroup} достаточно показать, что $H \prn HV$, а в силу леммы \ref{CrterionOfProSub} для этого достаточно, чтобы выполнялось равенство $U=N_U(H)[H,U]$ для любой $H$-инвариантной подгруппы $U\le V$.

Пусть $U$ --- $H$-инвариантная подгруппа из $V$. Включение $N_U(H)[H,U]\subseteq U$ очевидно. Заметим, что $S \le B$ является транзитивной группой подстановок на множестве $\{V_1, \ldots, V_n\}$. Ввиду леммы \ref{AwrH} имеем $U=C_U(S)[S,U]$. Теперь ясно, что $[H,U] \ge [S,U]$. Покажем, что $C_U(H) \ge C_U(S)$.
Это так, поскольку из леммы \ref{TransitiveAction} следует, что
$C_V(S) \le C_V(H)$ и $$C_U(S)=U\cap C_V(S) \le U\cap C_V(H)=C_U(H).$$

Теперь легко понять, что $$U=C_U(S)[S,U]\le C_U(H)[H,U]\le N_U(H)[H,U].$$

Значит, $H \prn HV$.  \qed\medskip

\begin{Lemma}\label{A4wrSn} Пусть $G=A\wr B$~--- естественное подстановочное сплетение группы $A \cong A_4$ и симметрической группы $B=S_n$, где $n=2^w$, и $H$ --- подгруппа нечетного индекса в $G$. Тогда $H \prn G$.
\end{Lemma}

\noindent{\sc Доказательство.} Заметим, что $G/O_2(G) \cong G_1$, где $G_1$ --- естественное подстановочное сплетение циклической группы $A$ порядка $3$ и симметрической группы $S_n$. Применение лемм~\ref{Quot} и \ref{3wrSn} завершает доказательство леммы. \qed\medskip

\begin{Lemma}\label{Zavarn} Пусть $Q$ --- подгруппа нечетного индекса группы $L=L_1\times L_2 \times \ldots \times L_n$, где $L_i$ --- конечные группы, и $\pi_i$ --- отображение проекции группы $L$ на $L_i$. Тогда если для некоторого $i$ группа $L_i$ проста и $Q^{\pi_i}=L_i$, то $L_i \le Q$.
\end{Lemma}

\noindent{\sc Доказательство.} Так как $L_i\unlhd L$, заключаем, что $Q\cap L_i\unlhd Q$, откуда $(Q\cap L_i)^{\pi_i}$ является нормальной подгруппой в группе $Q^{\pi_i}=L_i$, которая проста. Таким образом, либо $(Q\cap L_i)^{\pi_i}=1$, либо $(Q\cap L_i)^{\pi_i}=L_i$. В последнем случае получаем, что $L_i \le Q$.

Осталось показать, что $(Q\cap L_i)^{\pi_i}\not =1$. Пусть подгруппа $S \in \Syl_2(L)$ такова, что $S \le Q$. Тогда $S \cap L_i \in \Syl_2(L_i)$, откуда $1 \not =  S\cap L_i =(S \cap L_i)^{\pi_i} \le (Q \cap L_i)^{\pi_i}$. \qed\medskip

\begin{Lemma}\label{Conctruct}
Пусть $G=A\wr S_n=LH$~--- естественное подстановочное сплетение неабелевой простой группы $A$ и симметрической группы $H=S_n$, где $n=2^w$ и через $L=L_1\times \ldots \times L_n$, где $L_i \cong A$, обозначена база сплетения. Пусть $K$ --- подгруппа нечетного индекса в $G$ и $K_0=K \cap L$. Тогда для любой подгруппы $M_1$ группы $L_1$, которая содержит нормализатор в $L_1$ проекции группы $K_0$ на $L_1$, подгруппа $K$ содержится в некоторой подгруппе группы $G$, изоморфной $M_1 \wr  S_n$.

\end{Lemma}

\noindent{\sc Доказательство.} Группа $H=S_n$ естественным образом действует на множестве $\Omega=\{1, \ldots, n\}$.
Обозначим через $X_i$ стабилизатор точки $i$ в $H$.

Пусть $S$ --- силовская $2$-подгруппа группы $G$, содержащаяся в $K$. Тогда, поскольку все силовские $2$-подгруппы сопряжены в $G$,
существует элемент $g \in G$ такой, что подгруппа $S^{g}$ содержит некоторую регулярную подгруппу  $T$ группы $H=S_n$.
Поэтому мы можем считать, что $T \le K$.  Тогда $H=TX_i=X_iT$ для любого $i$.

Заметим, что %, поскольку $n=2^w$,
$T$ действует транзитивно на множестве $\{L_i\mid 1 \le i \le n\}$.
Пусть $K_i$~--- проекция группы $K_0$ на~$L_i$. Тогда $T$ также действует транзитивно на множестве $\{K_i \mid 1 \le i \le n\}$ (откуда $K_i \cong K_j$ для любых $i$ и $j$), поэтому \begin{equation}\label{eq1} \{K_i\mid 1 \le i \le n\}=\{K_1^t \mid t \in T\} \end{equation} и $K_0 \le \prod_{i=1}^n K_i$. Более того,  $K \le N_G\left(\prod_{i=1}^n K_i\right)$.

Заметим, что $N_L\left(\prod_{i=1}^n K_i\right)=\prod_{i=1}^n N_{L_i}(K_i)$. Покажем, что $H \le N_G\left(\prod_{i=1}^n K_i\right)$. Поскольку подгруппа $X_i$ централизует $L_i$, она централизует также $K_i$. Кроме того, $T$ регулярно действует на $\{K_i\mid 1 \le i \le n\}$
ввиду (\ref{eq1}). Для любого элемента $h \in H$ и для любого $i \in \{1, \ldots, n\}$ найдутся элементы $x \in X_i$ и $t \in T$ такие, что $h=xt$. Пусть $K_i^t=K_j$. Тогда $$K_i^h=K_i^{xt}=K_i^t=K_j.$$
Тем самым показано, что $H \le N_G\left(\prod_{i=1}^n K_i\right)$. Поэтому $$N_G\big(\prod_{i=1}^n K_i\big)=HN_L\big(\prod_{i=1}^n K_i\big) =H\prod_{i=1}^n N_{L_i}(K_i)\cong N_{L_1}(K_1)\wr H.$$

Кроме того, если $N_{L_1}(K_1)\le M_1 \le L_1$, то $$K \le H\prod_{i=1}^n N_{L_i}(K_i)=H\prod_{t \in T}\left(N_{L_1}(K_1)\right)^t \le H\prod_{t \in T}(M_1)^t \cong M_1 \wr H.$$

\qed\medskip

\begin{Lemma}\label{A5wrSn} Пусть $G=A\wr B$~--- естественное подстановочное сплетение группы $A=A_5$ и симметрической группы $B=S_n$, где $n=2^w$, и $H$ --- подгруппа нечетного индекса в $G$. Тогда $H \prn G$.
\end{Lemma}

\noindent{\sc Доказательство.} Обозначим базу сплетения через $L=L_1\times L_2 \times \ldots \times L_n$, где $L_i \cong A$. Пусть $S \in \Syl_2(G)$ такая, что $S \le H$. Заметим, что $G=X/O_2(X)$, где
$$
 X=SL_2(5)\wr S_{n}\hookrightarrow G_1=Sp_{2n}(5).
$$
Из \cite[теорема~9]{Maslova} следует, что индекс $|{G_1:X}|$ нечетен. Ввиду леммы \ref{Norm_Syl_Simple} выполнено равенство $|{N_{G_1}(S_1):S_1}|=3$, где $S_1 \in\Syl_2(G_1)$. Поэтому $|N_X(S_1):S_1|=|G:N_G(S)|$ делит $3$. Заметим, что элемент $(g,g,\ldots,g) \in L$, где $g$ --- элемент порядка $3$ из нормализатора в группе $A$ ее силовской $2$-подгруппы,  принадлежит $N_G(S)$. Поэтому $|N_G(S):S|=3$.

Пусть $g$ --- элемент порядка $3$ из $N_G(S)$. Приведенные рассуждения показывают, что $g \in L$.  Покажем, что подгруппы $H$ и $H^g$ сопряжены в $K=\langle H, H^g\rangle$.

Рассмотрим подгруппу $K_0=K\cap L$. Пусть $\pi_i$ --- отображение координатной проекции группы $L$ на $L_i$. Заметим, что $K$ действует транзитивно на $L_i$, поэтому все подгруппы $K_0^{\pi_i}$ попарно изоморфны. Кроме того, $K_0^{\pi_i}$ являются подгруппами нечетных индексов в $L_i$,
поэтому для них есть всего три возможности:

$(1)$ $K_0^{\pi_i} \cong A_5$ для любого $i${\rm;}

$(2)$ $K_0^{\pi_i} \cong A_4$ для любого $i${\rm;}

$(3)$ $K_0^{\pi_i} \cong C_2 \times C_2$ для любого $i$.

В случае $(1)$ ввиду леммы \ref{Zavarn} подгруппа $L$ содержится в $K$, поэтому $g \in K$.

В случае $(2)$ подгруппа $K$ нормализует подгруппу $\prod_{i=1}^n K_0^{\pi_i}$ и $N_G(S) \le N_G(\prod_{i=1}^n K_0^{\pi_i})$.
Поэтому ввиду леммы \ref{Overgroup} достаточно показать, что $H \prn N_G(\prod_{i=1}^n K_0^{\pi_i})$.
Это так, потому что ввиду леммы \ref{A4wrSn} подгруппы нечетных индексов пронормальны в группе $N_G(\prod_{i=1}^n K_0^{\pi_i}) \cong A_4\wr S_n$.

В случае $(3)$ имеем $K\cap L=H \cap L=S \cap L \in \Syl_2(L)$ и $L\cap S \unlhd H$. Поэтому $H \le N_G(L \cap S)$ и $g \in N_G(L \cap S)$ ввиду выбора элемента $g$. Заметим, что подгруппы нечетных индексов пронормальны в группе $N_G(L \cap S)/(L \cap S) \le B \cong S_n$. Следовательно, подгруппа $H$ и $H^g$ сопряжены  в $\langle H, H^g \rangle$ ввиду леммы \ref{Quot}. Поэтому $H \prn G$. \qed\medskip

\begin{Lemma}\label{PSL2wrSn} Пусть $G=A\wr B$~--- естественное подстановочное сплетение группы $A=PSL_2(q)$ и симметрической группы $B=S_n$, где $q \equiv \pm 3 \pmod 8$ и $n=2^w$, и ${H}$ --- подгруппа нечетного индекса в $G$. Тогда ${H} \prn G$.
\end{Lemma}

\noindent{\sc Доказательство.} Допустим, что лемма неверна и $q$~--- наименьшее из чисел, сравнимых с $\pm3$ по модулю 8,
для которых группа $G$ содержит непронормальную подгруппу ${H}$ нечетного индекса. Пусть подгруппа $S \in \Syl_2(G)$ такая, что $S \le {H}$.
%Можно считать, что $S$ содержит некоторый цикл $n_H$ длины $n$ из $B$.

Обозначим базу сплетения через $L=L_1\times L_2 \times \ldots \times L_n$, где $L_i \cong A$.

Заметим, что $G=X/O_2(X)$, где
$$
 X=SL_2(q)\wr S_{n}\hookrightarrow G_1=Sp_{2n}(q).
$$
Из \cite[теорема~9]{Maslova} следует, что индекс $|G_1:X|$ нечетен. Теми же рассуждениями, что и в доказательстве леммы~\ref{A5wrSn} убеждаемся, что $|N_G(S):S|=3$.

Пусть $g$ --- элемент порядка $3$ из $N_G(S)$. Ввиду сказанного,  $g \in L$. Более того, можно считать, что $g \in C_L(B)$. Покажем, что подгруппы ${H}$ и ${H}^g$ сопряжены в подгруппе $K=\langle {H}, {H}^g\rangle$ и тем самым по лемме~\ref{NormSyl} докажем пронормальность в $G$ подгруппы~$H$.

Рассмотрим подгруппы $H_0=H\cap L$ и $K_0=K\cap L$. Пусть $\pi_i: L\rightarrow L_i$ --- отображение координатной проекции. Заметим, что так как $H$ содержит силовскую 2-подгруппу группы $G$,  $H$ действует сопряжениями транзитивно на множестве $\{L_i\mid i=1,\dots,n\}$. Поэтому все группы $H_i=H_0^{\pi_i}$ попарно изоморфны и аналогичное утверждение справедливо для групп $K_i=K_0^{\pi_i}$. Кроме того, $H_i$ и $K_i$ являются подгруппами нечетных индексов в $L_i$, поскольку содержат силовскую 2-подгруппу $S_i=(S\cap L)$ группы~$L_i$.

Существует две возможности.

\smallskip

С л  у ч а й~~$(i)$: $K_i = L_i$ для любого $i$.

\smallskip

С л  у ч а й~~$(ii)$: $K_i < L_i$ для любого $i$.

\smallskip

В случае $(i)$ ввиду леммы \ref{Zavarn} имеем $L \le K$, поэтому $g \in K$ и требуемое доказано.

Пусть имеет место случай $(ii)$. Ввиду леммы \ref{Conctruct}
можно считать, что $K \le M_1 \wr S_n$, где $K_1 \le M_1 < L_1$, и $M_1$ --- максимальная подгруппа нечетного индекса в $L_1$.

По \cite[теорема 1]{Maslova} имеет место один из следующих четырех случаев.

\medskip
С л  у ч а й~~$(ii)(1)$: $M_1\cong A_4$. В этом случае ${H} \prn M_1 \wr S_n$ ввиду леммы \ref{A4wrSn} и $N_G(S) \le M_1 \wr S_n$, поэтому ${H} \prn G$ ввиду леммы \ref{Overgroup}.

\medskip
С л  у ч а й~~$(ii)(2)$: $M_1\cong A_5$. В этом случае ${H} \prn M_1 \wr S_n$ ввиду леммы \ref{A5wrSn} и $N_G(S) \le M_1 \wr S_n$, поэтому ${H} \prn G$ ввиду леммы \ref{Overgroup}.

\medskip
С л  у ч а й~~$(ii)(3)$: $M_1\cong PSL_2(q_0)$, где $q=q_0^r$ для некоторого нечетного простого числа~$r$.
Легко видеть, что $q_0\equiv \pm 3\pmod 8$, откуда следует, что ${H} \prn M \wr S_n$ ввиду выбора числа $q$, кроме того, $N_G(S) \le M \wr S_n$ и, значит,  ${H} \prn G$ по лемме \ref{Overgroup}.

\medskip
С л  у ч а й~~$(ii)(4)$: $M_1$~--- диэдральная группа порядка $q-\varepsilon$, где число $\varepsilon=\pm1$ выбрано так, что
$q\equiv \varepsilon\pmod 4$. Мы покажем сначала, что в этом случае  $g\in N_L(H_0)$.

Заметим, что в рассматриваемом случае группа $M_1$ и любая ее подгруппа (в частности $H_1$ и $K_1$ и, как следствие, все группы  $H_i$ и $K_i$) обладают нормальными циклическими 2-дополнениями. Так как $$H_0\le \langle H_1,\dots,H_n\rangle\cong H_1\times\dots\times H_n,$$
подгруппа $H_0$ также обладает нормальным $2$-дополнением $U=O(H_0).$

Положим $g_i=g^{\pi_i}$ для всех $i=1,\dots,n$.
Пусть также $U_i=O(H_i)$ и  $V_i=O(K_i)$ для любого~$i$. Тогда  $H_i^{g_i}=(H^g)^{\pi_i}\le K_i$. Поэтому $$U_i^{g_i}=O(H_i)^{g_i}=O(H_i^{g_i})\le O(K_i)=V_i.$$ Но $U_i$~--- единственная подгруппа группы $V_i$, порядок которой равен $|U_i|$, ввиду цикличности~$V_i$. Поэтому $U_i^{g_i}=U_i$, и согласно лемме \ref{OddOrderSubgroupsPSL} элемент $g_i\in L_i\cong PSL_2(q)$ централизует $U_i$.

Далее, поскольку $O(H_0)=U\le\langle U_1,\dots,U_n\rangle$, имеем $$g\in \langle g_1,\dots, g_n\rangle\le \langle C_{L_1}(U_1),\dots,C_{L_n}(U_n)\rangle\le C_{L}(\langle U_1,\dots,U_n\rangle)\le C_L(U).$$ Отсюда, поскольку $U$~--- 2-дополнение, а $S\cap L$~--- силовская 2-подгруппа в $H_0$, получаем $$H_0^g=(S^g\cap L)U^g=(S\cap L)U=H_0$$ и, как и утверждалось, $g\in N_L(H_0)$.

Рассмотрим в $G$ подгруппу $Y=N_G(H_0)$. По доказанному $g\in Y$ и ясно, что $H\le Y$, откуда $$K=\langle {H}, {H}^g\rangle\le\langle {H}, g\rangle\le Y.$$ Далее, $Y$~--- подгруппа нечетного индекса в $G$ и проекция группы $Y\cap L=N_L(H_0)$ на каждый сомножитель $L_i$ строго меньше $L_i$ (иначе из леммы \ref{Zavarn} следует, что $L_i\le Y$ и $L\le Y$, ввиду транзитивности действия группы $Y$ на множестве $\{L_i\mid i=1,\dots,n\}$, откуда следовало бы, что подгруппа $L$ содержит разрешимую нетривиальную нормальную подгруппу $H_0$). По лемме \ref{Conctruct} имеем $Y\le M< G$, где $M\cong  M_1^*\wr S_n$, а $M_1^*$~--- максимальная подгруппа группы $L_1$, содержащая $Y^{\pi_1}$. Для группы $M_1^*$, как и для группы $M_1$ выше, имеет место один из случаев $(ii)(1)$--$(ii)(4)$. Леммы \ref{DwrSn}, \ref{A4wrSn}, \ref{A5wrSn} и выбор $q$ показывают, что во всех случаях $H\prn M$ и, поскольку $g\in Y\le M$, подгруппы $H$ и $H^g$ сопряжены в $\langle {H}, {H}^g\rangle$.
\qed\medskip

\section*{Доказательство теоремы 1}

Пусть $G=Sp_{2n}(q)$, где  $q \equiv \pm 3 \pmod 8$ и $n$ не является числом вида $2^w$ или $2^w(2^{2k}+1)$.
Тогда  в двоичной записи числа $n$ либо найдутся  две единицы в некоторых разрядах c номерами $s_1$ и $s_2$ разной четности, либо три единицы в разрядах с номерами $s_1,$ $s_2$ и $s_3$  одинаковой четности.
Пусть $m=2^{s_1}+2^{s_2}$ или $m=2^{s_1}+2^{s_2}+2^{s_3}$ соответственно. Легко видеть, что $m$ кратно 3.

Рассмотрим стабилизатор $H$ в $G$ невырожденного подпространства
размерности $m$ пространства $V$. Индекс $|G:H|$ нечетен ввиду \cite[теорема~1]{Maslova}.  Таким образом, если $K$~--- подгруппа нечетного индекса в $H$, то $K$~--- подгруппа нечетного индекса в группе~$G$.

Легко понять, что $H = H_1 \times H_2$, где $H_1\cong Sp_{2m}(q)$ и $H_2\cong Sp_{2(n-m)}(q)$.
Так как $3$ делит $m$, подгруппа $H_1$ содержит по лемме~\ref{NonProNorm} непронормальную подгруппу $K_1$ нечетного индекса. Поэтому ввиду леммы~\ref{Quot} подгруппа $K_1\times H_2$ непронормальна в $H$, а значит, по лемме \ref{Overgroup} и в~$G$.

Ввиду леммы \ref{Quot} простая группа $G/O_2(G)=G/Z(G)\cong PSp_{2n}(q)$ содержит
непронормальную подгруппу нечетного индекса $(K_1 \times H_2)/O_2(G)$.
\qed\medskip

\section*{Доказательство теоремы 2}
Пусть $G = Sp_n(q)$, где $n=2^w$ и $q \equiv \pm 3 \pmod 8$.

Допустим, что теорема неверна и $q$~--- наименьшее из чисел, сравнимых с $\pm3$ по модулю 8,
для которых группа $G$ содержит непронормальную подгруппу $H$ нечетного индекса.

Рассмотрим естественный гомоморфизм $$\bar \,: G \rightarrow G/Z(G).$$ Заметим, что $|Z(G)|=2$, поэтому ввиду
леммы \ref{Quot} и выбора $H$ подгруппа $\bar H$ непронормальна $\bar G$. Поскольку ввиду
\cite[теорема]{KMR_1} подгруппы нечетных индексов пронормальны в группах $PSp_2(q) \cong PSL_2(q)$,
получаем $n \ge 4$.

Пусть $S \in \Syl_2(H)\subseteq \Syl_2(G)$. Ввиду леммы \ref{Norm_Syl_Simple} имеем $|N_{\bar G}(\bar S):\bar S|=3$.
По лемме~\ref{NormSyl} найдется  элемент $\bar g$ порядка $3$ из $N_{\bar G}(\bar S) \setminus \bar S$  такой, что $\bar H$ и $\bar H^{\bar g}$ не сопряжены в $\bar K =\langle \bar H, \bar H^{\bar g}
 \rangle$.

В силу выбора элемента $\bar g$ подгруппа $\bar K$ собственная в $\bar G$. Значит, существует максимальная подгруппа $\bar M$ (нечетного индекса) в $\bar G$ такая, что $\bar K \le \bar M$.

Ввиду \cite[теорема~1]{Maslova} для $\bar M$ существуют следующие возможности.

\medskip

С л  у ч а й~~$(1)$: $\bar M = C_{\bar G}(\sigma)$ для полевого автоморфизма $\sigma$ простого нечетного порядка $r$ группы $\bar G$.

В этом случае ввиду \cite[предложение~4.5.4]{KL} имеем $M \cong PSp_n(q_0)$, где $q=q_0^r$. Поскольку $r$~--- нечетное число,
легко понять, что $q_0 \equiv \pm 3 \pmod 8$. Отсюда следует, что $\bar H \prn \bar M$ ввиду выбора числа $q$. Кроме того, $\bar S \le \bar M$ и $|N_{\bar M}(\bar S):\bar S|=3$ ввиду леммы~\ref{Norm_Syl_Simple}. Поэтому $\bar g \in N_{\bar G}(\bar S)=N_{\bar M}(\bar S) \le \bar M$.  Получаем противоречие между выбором $\bar g$ и тем, что $\bar H \prn \bar M$.

\medskip

С л  у ч а й~~$(2)$: $\bar G = PSp_4(q)$ и $\bar M \cong 2^4.A_5$, где  $q \equiv 3 \pmod 8$. Из \cite[теорема]{KMR_1} следует, что подгруппы нечетных индексов пронормальны в $\bar M/O_2(\bar M)$.  Отсюда  $\bar H \prn \bar M$ ввиду леммы~\ref{Quot}. Кроме того, $\bar S \le \bar M$ и $|N_{\bar M}(\bar S):\bar S|=3$ ввиду леммы~\ref{Norm_Syl_Simple}. Снова $\bar g \in N_{\bar G}(\bar S)=N_{\bar M}(\bar S) \le \bar M$, и снова противоречие между выбором $\bar g$ и тем, что $\bar H \prn \bar M$.

\medskip

С л  у ч а й~~$(3)$: $\bar M$ --- стабилизатор в $\bar G$ ортогонального разложения
\begin{equation}\label{bigoplus}
V = \bigoplus\limits_{i=1}^t V_i
\end{equation}
  естественного симплектического модуля $V$ группы $G$ в прямую
сумму изометричных подпространств $V_i$ размерности $m$, причем $m = 2^{w_1} \ge 2$. Выберем подгруппу  $\bar M$ указанного типа так, чтобы число $m$ было наименьшим возможным, т.~е. так, чтобы разложение~(\ref{bigoplus}) было неизмельчаемым.

Ввиду \cite[предложение~4.2.10]{KL} имеем $$\bar M \cong 2^{t-1}.(PSp_m(q) \wr S_t), \mbox{ где } n=mt.$$
Кроме того, $\bar S \le \bar M$ и в $\bar M$ существует элемент порядка $3$, нормализующий подгруппу $\bar S$, т.\,е. $|N_{\bar M}(\bar S):\bar S|=3$ ввиду леммы~\ref{Norm_Syl_Simple}. Поэтому, как и в предыдущих случаях, $\bar g \in N_{\bar G}(\bar S)=N_{\bar M}(\bar S) \le \bar M$.

 Рассмотрим композицию $\tilde{\,\,}:  M \rightarrow \bar M/O_2(\bar M)$ естественных гомоморфизмов
 $$
 M \rightarrow \bar M\,\,\,\text{ и }\,\,\,\bar M\rightarrow \bar M/O_2(\bar M).
 $$
  %Ввиду леммы~\ref{Quot} имеем $\bar H \prn \bar M$ тогда и только тогда, когда $\tilde{H} \prn \tilde{M};
  Строение группы $\bar M$ показывает, что $\tilde{M} \cong PSp_m(q) \wr S_t$.
  %По теореме Силова, не уменьшая общности, можно считать, что $\tilde{H}$ содержит цикл длины $t$ из $S_t$.

Обозначим через $L=L_1\times L_2 \times \ldots \times L_t$, где $L_i \cong PSp_m(q)$, нормальную  подгруппу в $\tilde{M}$, соответствующую базе сплетения $PSp_m(q) \wr S_t$. Заметим, что невырожденные подпространства $V_i$  в разложении~(\ref{bigoplus}) можно считать ассоциированными с соответствующими подгруппами $L_i$, рассматриваемыми как проективные симплектические группы.

Ввиду выбора элемента $\bar g$ его образ $\tilde g$ --- элемент порядка $3$ из $N_{\tilde{M}}(\tilde{S})$. Ясно, что $\tilde g \in L$. Покажем, что подгруппы $\tilde{H}$ и $\tilde{H}^{\tilde g}$ сопряжены в подгруппе $\tilde{K}=\langle \tilde{H}, \tilde{H}^{\tilde g}\rangle$. Тем самым мы получим противоречие, поскольку ядро гомоморфизма $\tilde{\,\,}$ содержится в $H$ и поэтому подгруппы $\bar H$ и $\bar H^{\bar g}$ окажутся сопряженными в $\bar K$.

Рассмотрим подгруппу ${K_0}=\tilde{K}\cap L$. Пусть $\pi_i$ --- отображение координатной проекции группы $L$ на $L_i$. Заметим, что $\tilde K$ действует транзитивно на $L_i$. Поэтому все $K_0^{\pi_i}$ попарно изоморфны. Кроме того, $K_0^{\pi_i}$ являются подгруппами нечетных индексов в~$L_i$.
Существуют две возможности.

\medskip

С л  у ч а й~~$(3)(i)$: $K_0^{\pi_i} = L_i$ для любого $i$.

\medskip

С л  у ч а й~~$(3)(ii)$: $K_0^{\pi_i} < L_i$ для любого $i$.

\medskip

Случай $(3)(i)$ исключается леммой \ref{Zavarn}, поскольку иначе $L \le \tilde{K}$ и $\tilde g \in \tilde{K}$.

\medskip

Пусть имеет место случай $(3)(ii)$.

Ввиду леммы \ref{Conctruct} можно считать, что $\tilde{K} \le Y
< \tilde M$, где  $Y\cong R_1 \wr S_t$, а  $R_1$ --- максимальная подгруппа нечетного индекса в $L_1\cong PSp_m(q)$.

Из леммы \ref{PSL2wrSn} следует, что $m>2$. Кроме того, поскольку подгруппа $\bar M$, содержащая $\bar K$,  была выбрана таким образом, чтобы число $m$ было наименьшим возможным, подгруппа $R_1$ в $L_1$ не может являться стабилизатором разложения ассоциированного с $L_1$ подпространства $V_1$  в ортогональную прямую сумму подпространств меньшей размерности. Поэтому согласно \cite[теорема 1]{Maslova} имеет место один из следующих  случаев.

\medskip
С л  у ч а й~~$(3)(ii)(1)$: $R_1 = C_{L_1}(\sigma)$ для полевого автоморфизма $\sigma$ простого нечетного порядка $r$ группы $L_1$, отождествляемой с $PSp_m(q)$.
Тогда ввиду \cite[предложение~4.5.4]{KL} имеем $R_1 \cong PSp_m(q_0)$, где $q=q_0^r$. Поскольку $r$ --- нечетное число,
легко понять, что $q_0 \equiv \pm 3 \pmod 8$. Отсюда следует, что $Y=X/O_2(X)$, где
$$ X=Sp_m(q_0)\wr S_{t}\hookrightarrow G_1=Sp_{m}(q_0).$$
Подгруппы нечетных индексов пронормальны в $Sp_{m}(q_0)$ ввиду выбора числа $q$.
Из \cite[теорема~9]{Maslova} следует, что индекс $|G_1:X|$ нечетен. Ввиду леммы \ref{Norm_Syl_Simple} имеем $|N_{G_1}(S_1):S_1|=3$, где $S_1 \in\Syl_2(G_1)$. Поэтому подгруппы нечетных индексов пронормальны в $X$ ввиду леммы \ref{Overgroup}, следовательно подгруппы нечетных индексов пронормальны в $Y$ ввиду леммы \ref{Quot}. Заметим, что подгруппа $Y
\cong R_1 \wr S_t$ содержит элемент порядка $3$, нормализующий ее силовскую $2$-подгруппу. Отсюда, как и выше, заключаем, что $\tilde g \in N_{\tilde{M}}(\tilde{S})\le Y$ и, так как $\tilde H\le \tilde K\le Y$, снова получаем противоречие.

\medskip
С л  у ч а й~~$(3)(ii)(2)$: $m=4$, $q$ --- простое число, $q \equiv 3 \pmod 8$ и $ R_1 \cong 2^4.A_5$.
Из леммы \ref{A5wrSn} следует, что подгруппы нечетных индексов пронормальны в $Y/O_2(Y)$.  Отсюда следует, что $\tilde{H} \prn Y$ ввиду леммы~\ref{Quot}. Кроме того, $\tilde{S} \le Y$ и, поскольку $Y$ содержит элемент порядка $3$, нормализующий $\tilde{S}$, имеем $\tilde g\in N_{\tilde{M}}(\tilde{S}) \le Y$.   Получаем противоречие.
\qed\bigskip

Авторы выражают глубочайшую благодарность А.~В.~Заварницину за полезные консультации.

%%%%%%%%%%%%%%%%%%% Список литературы %%%%%%%%%%%%%%%%%%%%%%%%%%%%

%%%%%%%%%%%%%%%%%%% Адрес автора %%%%%%%%%%%%%%%%%%%%

 Адреса авторов:

\smallskip

Кондратьев Анатолий Семёнович, Маслова Наталья Владимировна

Институт математики и механики им. Н. Н. Красовского УрО РАН,

ул. С. Ковалевской, 16,  Екатеринбург, 620990,

Уральский федеральный университет,

ул. Мира, 19, Екатеринбург, 620002,

a.s.kondratiev@imm.uran.ru, butterson@mail.ru

\medskip

Ревин Данила Олегович

Институт математики им. С.Л.Соболева СО РАН,

пр. акад. Коптюга, 4, Новосибирск, 630090,

Новосибирский государственный университет,

ул. Пирогова, 2, Новосибирск, 630090,

Department of Mathematics, University of Science and Technology of China,

Hefei 230026, P. R. China,

revin@math.nsc.ru

\vspace{3cm}


\begin{thebibliography}{100}

\addcontentsline{toc}{chapter}{Список
                  литературы}

{\small

\bibitem{Pronorm_Hall_Simple}
 {Е.~П.~Вдовин, Д.~О.~Ревин} Пронормальность холловых подгрупп в конечных простых группах // Сиб. матем. журн. 2012. Т~53, \No~3. С.~527--542.


\bibitem{KMR_1}
{А.~С.~Кондратьев, Н. В. Маслова,  Д.~О.~Ревин} О пронормальности подгрупп нечетного индекса в конечных простых группах //
Сиб. матем. журн. 2015. Т~56, \No~6. С.~1101--1107.


\bibitem{KMR_2}
{А.~С.~Кондратьев, Н. В. Маслова,  Д.~О.~Ревин} Критерий пронормальности добавлений к абелевым нормальным подгруппам //
Тр. Ин-та математики и механики УрО РАН. 2016. Т.~22, \No~1. C. 153--158.

%\bibitem{ProHall} {Д.~О.~Ревин, Е.~П.~Вдовин} О пронормальности холловых подгрупп // Сиб. матем. журн. 2013. Т.~54, \No~1. С.~35--43.

%\bibitem{Z} {G.~Glauberman} Central elements in core-free groups // J. Algebra. 1966. Vol.~4, \No~3. P.~403--420.

\bibitem{Atlas}
{J.~H.~Conway, R.~T.~Curtis, S.~P.~Norton, R.~A.~Parker, R.~A.~Wilson}
Atlas of finite groups. Oxford: Clarendon Press, 1985. 252~p.

\bibitem{KL}
{P.~B.~Kleidman,  M.~Liebeck} The subgroup structure of the finite classical groups. Cambridge: Cambridge University Press, 1990. 303~p.


\bibitem{Kond1}
{А.~С.~Кондратьев } Нормализаторы силовских 2-подгрупп в конечных простых группах //
 Мат. заметки. 2005. Т.~78, \No~3. С.~368--376.

\bibitem{Hupp} B. Huppert, Endliche Gruppen. Berlin: Springer-Verlag,  1967.

\bibitem{Isaacs}
{I.~M.~Isaacs} Character theory of finite groups. NY: Academic Press, 1976.

\bibitem{Isaacs1} {I.~M.~Isaacs} Finite Group Theory. Amer. Math. Soc.: Providence, Rhode Island, 2008.


\bibitem{Maz}
{В.~Д.~Мазуров } О множестве порядков элементов конечной группы  // Алгебра и логикаю. 1994. Т.~33, \No~1. C.~81--89.


%\bibitem{Rob} {D.~Robinson}, A course in the theory of groups. N. Y.: Springer-Verlag, 1996. 499~p.

%\bibitem{Hupp} {B.~Huppert} Endliche Gruppen I. Berlin: Springer, 1967. 794~p.

%\bibitem{KM} {А.~С.~Кондратьев, В.~Д.~Мазуров} 2-сигнализаторы конечных простых групп // Алгебра и логика. 2003. Т.~42, \No~5. С.~594--623.

%\bibitem{GLS} {D.~Gorenstein, R.~Lyons, R.~Solomon} The classification of the finite simple groups, Number~3 // Mathematical Surveys and Monographs. 1994. Vol.~40, \No~3. 419~p.

\bibitem{Maslova}
{Н. В. Маслова} Классификация максимальных подгрупп нечетного индекса в конечных простых классических группах //
Тр. Ин-та математики и механики УрО РАН. 2008. Т.~14, \No~4. C. 100--118.



\bibitem{CF} {R. Carter, P. Fong}, The Sylow 2-subgroups of the finite classical groups // J. Algebra. 1964. Т.~1, \No~1. C.~139--151.


}

\medskip
\end{thebibliography}
\end{document}